\documentclass[12pt,twoside]{article}
\usepackage{graphicx}
\usepackage{amsfonts}
\usepackage{amsbsy}
\usepackage{amssymb}
\usepackage{amsmath}
\usepackage{cite}
\usepackage{pst-all}
\usepackage{pstricks}
\usepackage{graphics}
\usepackage{tikz-cd}

\def\bpsp{\begin{pspicture}}
\def\epsp{\end{pspicture}}

\newtheorem{theorem}{Theorem}[section]
\newtheorem{remark}[theorem]{Remark}
\newtheorem{example}[theorem]{Example}
\newtheorem{lemma}[theorem]{Lemma}
\newtheorem{corollary}[theorem]{Corollary}
\newtheorem{definition}[theorem]{Definition}
\newtheorem{proposition}[theorem]{Proposition}
\newtheorem{note}{Note}
\newtheorem{case}{Case}
\newtheorem{conjecture}{Conjecture}
\newtheorem{question}{Question}

\newcommand{\bea}{\begin{eqnarray}}
\newcommand{\eea}{\end{eqnarray}}
\newcommand{\beq}{\begin{eqnarray*}}
\newcommand{\eeq}{\end{eqnarray*}}

\def\m4{\mbox{\rm ~(mod $4$)}}

\def \bd{\begin{definition}}
\def \ed{\end{definition}}
\def \be{\begin{example}}
\def \bp{\begin{proof}}
\def \ep{\end{proof}}
\def \ee{\end{example}}
\def \ed{\end{definition}}
\def \ed{\end{definition}}
\def \bqu{\begin{question}}
\def \equ{\end{question}}
\def \bcc{\begin{conjecture}}
\def \ecc{\end{conjecture}}
\def \bt{\begin{theorem}}
\def \et{\end{theorem}}
\def \bl{\begin{lemma}}
\def \el{\end{lemma}}
\def \bc{\begin{corollary}}
\def \ec{\end{corollary}}
\def \be{\begin{equation}}
\def \ee{\end{equation}}
\def \ben{\begin{enumerate}}
\def \een{\end{enumerate}}
\def \ba{\begin{array}}
\def \ea{\end{array}}
\def \bp{\begin{proposition}}
\def \ep{\end{proposition}}
\def \bx{\begin{example}}
\def \ex{\end{example}}
\def \br{\begin{remark}}
\def \er{\end{remark}}
\def \bdsc{\begin{description}}
\def \edsc{\end{description}}

\def \bn{\begin{case}}
\def \en{\end{case}}
\def \bnt{\begin{note}}
\def \ent{\end{note}}
\def\1{1\!\!1}

\def\mm2{\mbox{\rm ~(mod $2$)}}
\def\m4{\mbox{\rm ~(mod $4$)}}

\def\qed{\nolinebreak\hfill\rule{.2cm}{.2cm}\par\addvspace{.5cm}}

\def\m{\mu}

\def\1{\textbf{1}}
\def\0{\textbf{0}}

\begin{document}
\title{Total perfect codes in graphs realized by commutative rings}
\author{Rameez Raja $^{a, b}$\\
{\em $^{a}$ Department of Mathematics, National Institute of Technology,}\\ {Hazratbal, Srinagar 190006, India}\\
{\em $^{b}$ School of Mathematics, Harish-Chandra Research Institute, HBNI,}\\ {Chhatnag Road, Jhunsi, Allahabad 211019, India}\\
rameeznaqash@gmail.com, rameeznaqash@nitsri.ac.in}
\date{}

\pagestyle{myheadings} \markboth{Rameez}{Total perfect codes in graphs realized by commutative rings}

\maketitle

\noindent{\footnotesize \bf Abstract.}
Let $R$ be a commutative ring with unity not equal to zero and let $\Gamma(R)$ be a zero-divisor graph realized by $R$. For a simple, undirected, connected graph $G = (V, E)$, a {\it total perfect code} denoted by $C(G)$ in $G$ is a subset $C(G) \subseteq V(G)$ such that $|N(v) \cap C(G)| = 1$ for all $v \in V(G)$, where $N(v)$ denotes the open neighbourhood of a vertex $v$ in $G$. In this paper, we study total perfect codes in graphs which are realized as zero-divisor graphs. We show a zero-divisor graph realized by a local commutative ring with unity admits a total perfect code if and only if the graph has degree one vertices. We also show that if $\Gamma(R)$ is a regular graph on $|Z^*(R)|$ vertices, then $R$ is a reduced ring and $|Z^*(R)| \equiv 0(mod ~2)$, where $Z^*(R)$ is a set of non-zero zero-divisors of $R$. We provide a characterization for all commutative rings with unity of which the realized zero-divisor graphs admit total perfect codes. Finally, we determine the cardinality of a total perfect code in $\Gamma(R)$ and discuss the significance of the study of total perfect codes in graphs realized by commutative rings with unity.

\vskip 3mm



\section{\bf Introduction}

The notions of perfect codes in graphs evolved from the work \cite{B}, which in turn has roots in coding theory \cite{VL}. The theory of perfect codes forms an interesting part of combinatorics and has connections in group theory, diophantine number theory, geometry and cryptography. Perfect codes played a central role in the fast growing of error-correcting codes theory. Hamming and Golay \cite{G, H} constructed perfect binary single-error correcting codes of length $n$, where $n = 2k - 1$ for some integer $k$. Infinite classes of graphs with perfect codes have been constructed by Cameron, Thas and Payne \cite{CPT}, Thas \cite{T}, Hammond \cite{Ha} and the others. The existence of perfect codes has also been proved in Towers of Hanoi graphs \cite{CN}, Sierpinski graphs \cite{KMP} and Cayley graphs \cite{SZ}. For more on perfect codes in graphs, see \cite{AD, CRC, DS, GHT, LS, S}.\\

The study of combinatorial structures (graphs) realized by algebraic structures (rings) was first introduced by Beck \cite{Bk} and was further studied by Anderson and Livingston in \cite{AL}. The zero-divisor graph introduced in \cite{AL} is denoted by $\Gamma(R)$ with vertex set $Z^*(R) = Z(R)\setminus \{0\}$, and two vertices $x, y\in Z^*(R)$ are adjacent in $\Gamma(R)$ if and only if $xy = 0$. They believed that the study of combinatorial properties of zero-divisors better illustrates the zero-divisor structure of a ring. The zero-divisor graph of a commutative ring has also been studied in \cite{ALS, SR, SR1, SRR1} and was extended by Redmond \cite{Rd} to non-commutative rings. Redmond \cite{Rd*} also extended the zero-divisor graph of a commutative ring to an ideal-based zero-divisor graph of a commutative ring. The combinatorial properties of zero-divisors discovered in \cite{Bk} which have also been studied in module theory. Recently in \cite{SR}, the elements of a module $M$ have been classified into full-annihilators, semi-annihilators and star-annihilators. This classification of elements in $M$ has unfolded a correspondence between ideals in $R$, submodules of $M$ and vertices of annihilating graphs arising from $M$.

In a graph $G$, the distance between two vertices $u, v\in V(G)$ denoted by $d(u, v)$ is the length of a shortest path between $u$ and $v$ in $G$ and is defined to be infinity if no path between $u$ and $v$ exists. The \textit{ball} with centre $v\in V(G)$ and radius $t$ is the set of vertices of $G$ with distance atmost $t$ to $v$. A \textit{code} denoted by $C(G)$ in $G$ is simply a subset of $V(G)$. For an integer $t\geq 1$, a code $C(G) \subseteq V (G)$ is called a \textit{perfect $t$-code} \cite{K} in $G$ if balls with centres in $C(G)$ and radius $t$ form a partition of $V (G)$. A code $C(G)$ is said to be a \textit{total perfect code} \cite{GHT} in $G$ if every vertex of $G$ has exactly one neighbour in $C(G)$, that is, $|N(v) \cap C(G)| = 1$ for all $v \in V(G)$. A \textit{perfect 1-code} in a graph is also called an \textit{efficient dominating set} \cite{DS, KP} or \textit{independent perfect dominating set} \cite{L}, and a total perfect code is called an \textit{efficient open dominating set} \cite{HHS}. Total perfect codes in graphs have potential applications in some practical domains, such as placement of Input/Output devices in a supercomputing network so that each element to be processed is at a distance  of atmost one to exactly one Input/Output device \cite{AD}. It is known \cite{GSS} that deciding whether a graph has a total perfect code is NP-complete.\\ 

In Section 2, we discuss total perfect codes in graphs realized as zero-divisor graphs. From \cite{AL}, it can be explored and confirmed that the classes of graphs like paths, cycles, regular graphs, bipartite graphs, complete bipartite graphs etc. are realized as zero-divisor graphs. For a local commutative ring $R$ with unity, Propositions 2.3 and 2.6 are characterizations of graphs (realized by $R$) which admit total perfect codes. In Lemma 2.8, we show that if $\Gamma(R)$ is a regular graph, then $R$ is a reduced ring and $|Z^*(R)| \equiv 0(mod ~2)$. Lemma 2.10 shows that among all bipartite graphs, only complete bipartite graphs admit total perfect codes. Theorem 2.13 which is the main result of this paper provides a characterization for all commutative rings of which the realized zero-divisor graphs admit total perfect codes. Finally, we finish the paper with Theorem 2.17 in which we determine the cardinality of a total perfect code in $\Gamma(R)$.\\

We conclude this section with some notations, which are mostly standard and will be used throughout this paper.\\ 

Throughout, $R$ is a commutative ring (with $1 \neq 0$). We denote the ring of integers by $\mathbb{Z}$, non-negative integers by $\mathbb{Z}_{\geq0}$, positive integers by $\mathbb{N}$ and the ring of integers modulo $n$ by $\mathbb{Z}_n$.\\

\section{\bf Graphs realized by commutative rings}
In this Section, we discuss total perfect codes in graphs which are realized as zero-divisor graphs. Let $R$ realizes $\Gamma(R)$ as its zero-divisor graph with vertex set $Z^*(R)$. Then $\Gamma(R)$ is said to admit a total perfect code if there exits some set $C(R)\subseteq Z^{*}(R)$ such that $|Ann(x) \cap C(R)| = 1$ for all $x\in Z^*(R)$, where $Ann(x) = ann(x)\setminus \{0, x\}$ denotes the open neighbourhood of $x$ in $\Gamma(R)$ and $ann(x) = \{y\in R ~|~ xy = 0\}$ denotes the annihilator of an element $x$ in $R$. It follows that if $R$ realizes $\Gamma(R)$ as its zero-divisor graph, then a set $C(R)$ is a total perfect code in $\Gamma(R)$ if and only if the subgraph induced by $C(R)$ is a matching in $\Gamma(R)$ and the set $\{Ann(x) ~|~ x \in C(R) \}$ is a partition of $Z^*(R)$. In particular, any total perfect code in $\Gamma(R)$ contains an even number of vertices.\\

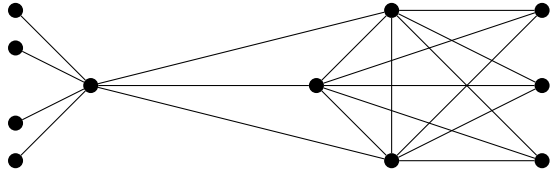
It is simple to construct zero-divisor graphs given a commuutative ring with unity. For example, let $R = \mathbb{Z}_2 \times \mathbb{Z}_2$. Then $Z^*(R) = \{(0, 1), (1, 0)\}$ and $\Gamma(R)$ is a path on two vertices. But the opposite endeavour is more complicated. That is, endow a commutative ring for the given graph. We see below that if $R$ realizes some graph as its zero-divsor graph, then it admits a trivial matching or empty set as its total perfect code. Thus, the study of total perfect codes in zero-divisor graphs associated with commutative rings simplifies the drudgery of determining a commutative ring for the given graph. It is worth to note that by [Theorem 2.3, \cite{AL}] the diameter of $\Gamma(R)$ is atmost $3$, so not every graph can be realized as a zero-divisor graph. However, there are well known families of graphs which are realized by some commutative rings as its zero-divisor graphs. $\mathbb{Z}_2[X_1, X_2, \cdots, X_n]/ (X^{2}_{1}, X^{2}_{2}, \cdots, X^{2}_{n}, X_1X_2, X_2X_3, \cdots, X_{n-1}X_{n}, X_nX_1)$ realizes a graph which contains a cycle on $X_1, X_2, \cdots, X_n$ vertices, $\mathbb{Z}_8$ realizes a path, $\mathbb{Z}_3 \times \mathbb{Z}_3$ realizes a regular graph, $\mathbb{Z}_2 \times \mathbb{Z}_5$ realizes a star graph, $\mathbb{Z}_2 \times \mathbb{Z}_4$ realizes a bipartite graph and $\mathbb{Z}_{25}$ realizes a complete graph as its zero-divisor graph.\\     

We begin the investigation of these codes with the following two examples.\\ 

\begin{example} Consider a graph $G_1$ as shown in Figure 1 with vertex set $V(G_1) = \{x_1, x_2, x_3, x_4,$
\noindent$ x_5, x_6, x_7\}$ and  edge set $E(G_1) = \{(x_1, x_2), (x_1, x_3), (x_1, x_4), (x_1, x_6), (x_4, x_5), (x_4, x_7), (x_5, x_6), (x_6, x_7)\}$. For a ring $\mathbb{Z}_{12}$, we have $Z^*(\mathbb{Z}_{12}) = \{\bar{2}, \bar{3}, \bar{4}, \bar{6}, \bar{8}, \bar{9}, \bar{10}\}$. Therefore $|V(G_1)| = |Z^*(\mathbb{Z}_{12})|$. $ann(\bar{6}) = \{ \bar{0}, \bar{2}, \bar{4}, \bar{6}, \bar{8}, \bar{10}\}$, whereas $Ann(\bar{6}) = \{ \bar{2}, \bar{4}, \bar{8}, \bar{10} \}$. Similalrly $Ann(\bar{2}) = \{ \bar{6} \}$, $Ann(\bar{3}) = \{ \bar{4}, \bar{8} \}$, $Ann(\bar{4}) = \{ \bar{3}, \bar{6}, \bar{8} \}$, $Ann(\bar{8}) = \{ \bar{3}, \bar{6}, \}$, $Ann(\bar{9}) = \{ \bar{4}, \bar{8}, \bar{10} \}$ and $Ann(\bar{10}) = \{\bar{6}\}$. It is clear from the adjacencies in $G_1$ that $Ann(\bar{6}) = N(x_1)$, $Ann(\bar{2}) = N(x_2)$, $Ann(\bar{3}) = N(x_5)$, $Ann(\bar{4}) = N(x_4)$, $Ann(\bar{8}) = N(x_6)$, $Ann(\bar{9}) = N(x_7)$ and $Ann(\bar{10}) = N(x_3)$. It follows from the above graph relations that $\mathbb{Z}_{12}$ realizes $G_1$ as its zero-divisor graph. It is easy to verify that there is a subset $C(\mathbb{Z}_{12}) = \{ \bar{6}, \bar{4}\}$ such that for all $x\in Z^{*}(\mathbb{Z}_{12})$ $| Ann(x) \cap C(\mathbb{Z}_{12})| = 1$. Therefore $\Gamma(\mathbb{Z}_{12})$ admits a total perfect code.
\end{example}
\begin{align*}
\begin{pgfpicture}{9cm}{-1cm}{3cm}{1.5cm}
\pgfnodecircle{Node1}[fill]{\pgfxy(7,0)}{0.1cm}
\pgfnodecircle{Node2}[fill]{\pgfxy(5.6,1)}{0.1cm}
\pgfnodecircle{Node3}[fill]{\pgfxy(5.6,-1)}{0.1cm}
\pgfnodecircle{Node4}[fill]{\pgfxy(10, 1)}{0.1cm}
\pgfnodecircle{Node5}[fill]{\pgfxy(9, 0)}{0.1cm}
\pgfnodecircle{Node6}[fill]{\pgfxy(10, -1)}{0.1cm}
\pgfnodecircle{Node7}[fill]{\pgfxy(11.2,0)}{0.1cm}
\pgfnodeconnline{Node1}{Node2}
\pgfnodeconnline{Node1}{Node3}
\pgfnodeconnline{Node1}{Node4}
\pgfnodeconnline{Node1}{Node6}
\pgfnodeconnline{Node5}{Node4}
\pgfnodeconnline{Node5}{Node6}
\pgfnodeconnline{Node6}{Node7}
\pgfnodeconnline{Node7}{Node4}
\pgfputat{\pgfxy(6.9,.3)}{\pgfbox[left,center]{$x_1$}}
\pgfputat{\pgfxy(5.5, 1.3)}{\pgfbox[left,center]{$x_2$}}
\pgfputat{\pgfxy(5.5, -1.4)}{\pgfbox[left,center]{$x_3$}}
\pgfputat{\pgfxy(9.8, 1.3)}{\pgfbox[left,center]{$x_4$}}
\pgfputat{\pgfxy(8.6, .3)}{\pgfbox[left,center]{$x_5$}}
\pgfputat{\pgfxy(9.8,-1.4)}{\pgfbox[left,center]{$x_6$}}
\pgfputat{\pgfxy(11.1, .3)}{\pgfbox[left,center]{$x_7$}} 
\end{pgfpicture}
\end{align*} \hskip .9cm \hskip .9cm \hskip .9cm \hskip .9cm \hskip .9cm \hskip .9cm \hskip .3cm Figure 1. $G_1 = \Gamma(\mathbb{Z}_{12})$.\\

\begin{example} Consider a graph $G_2$ as shown in Figure 2 with vertex set $V(G_2) = \{y_1, y_2, y_3, y_4,$
\noindent$ y_5, y_6, y_7, y_8, y_9, y_{10}, y_{11}\}$ and  edge set $E(G_2) = \{(y_1, y_2), (y_1, y_3), (y_1, y_4), (y_1, y_5), (y_1, y_6), (y_1, y_7),$ 
\noindent $(y_1, y_8), (y_8, y_9), (y_8, y_6), (y_8, y_7), (y_8, y_{10}), (y_8, y_{11}), (y_6, y_{11}), (y_7, y_{11}) \}$. For a ring $ \mathbb{Z}_2 \times \mathbb{Z}_8$, we have $Z^*(\mathbb{Z}_2 \times \mathbb{Z}_8) = \{(\bar{0}, \bar{1}), (\bar{0}, \bar{2}), (\bar{0}, \bar{3}), (\bar{0}, \bar{4}), (\bar{0}, \bar{5}), (\bar{0}, \bar{6}), (\bar{0}, \bar{7}), (\bar{1}, \bar{0}), (\bar{1}, \bar{2}), (\bar{1}, \bar{4}), (\bar{1}, \bar{6})\}$. Therefore $|V(G_2)| = |Z^*(\mathbb{Z}_2 \times \mathbb{Z}_8)|$ and $Ann(\bar{0}, \bar{6}) = N(y_1)$, $Ann(\bar{0}, \bar{1}) = N(y_2)$, $Ann(\bar{0}, \bar{3}) = N(y_3)$, $Ann(\bar{0}, \bar{5}) = N(y_4)$, $Ann(\bar{0}, \bar{7}) = N(y_5)$, $Ann(\bar{0}, \bar{2}) = N(y_6)$, $Ann(\bar{0}, \bar{6}) = N(y_7)$, $Ann(\bar{0}, \bar{4}) = N(y_8)$, $Ann(\bar{1}, \bar{2}) = N(y_8)$, $Ann(\bar{1}, \bar{6}) = N(y_{10})$, $Ann(\bar{1}, \bar{4}) = N(y_{11})$. Thus $\mathbb{Z}_2 \times \mathbb{Z}_8$ realizes $G_2$ as its zero-divisor graph. It can be verified from $\Gamma(\mathbb{Z}_2 \times \mathbb{Z}_8)$ that there is no subset $C(\mathbb{Z}_2 \times \mathbb{Z}_8)\subseteq Z^{*}(\mathbb{Z}_2 \times \mathbb{Z}_8)$ such that for all $x\in Z^{*}(\mathbb{Z}_2 \times \mathbb{Z}_8)$ $|Ann(x) \cap C(\mathbb{Z}_2 \times \mathbb{Z}_8)| = 1$.\\
\end{example}
\begin{align*}
\begin{pgfpicture}{9.3cm}{-1cm}{3cm}{0.9cm}
\pgfnodecircle{Node1}[fill]{\pgfxy(7,0)}{0.1cm}
\pgfnodecircle{Node2}[fill]{\pgfxy(6,1)}{0.1cm}
\pgfnodecircle{Node3}[fill]{\pgfxy(6, .5)}{0.1cm}
\pgfnodecircle{Node4}[fill]{\pgfxy(6, -.5)}{0.1cm}
\pgfnodecircle{Node5}[fill]{\pgfxy(6, -1)}{0.1cm}
\pgfnodecircle{Node6}[fill]{\pgfxy(10, 1)}{0.1cm}
\pgfnodecircle{Node7}[fill]{\pgfxy(10,-1)}{0.1cm}
\pgfnodecircle{Node8}[fill]{\pgfxy(9.5, 0)}{0.1cm}
\pgfnodecircle{Node9}[fill]{\pgfxy(10.5, 0.4)}{0.1cm}
\pgfnodecircle{Node10}[fill]{\pgfxy(10.5,-0.4)}{0.1cm}
\pgfnodecircle{Node11}[fill]{\pgfxy(12.3,0)}{0.1cm}
\pgfnodeconnline{Node1}{Node2}
\pgfnodeconnline{Node1}{Node3}
\pgfnodeconnline{Node1}{Node4}
\pgfnodeconnline{Node1}{Node5}
\pgfnodeconnline{Node1}{Node6}
\pgfnodeconnline{Node1}{Node7}
\pgfnodeconnline{Node1}{Node8}
\pgfnodeconnline{Node6}{Node8}
\pgfnodeconnline{Node8}{Node7}
\pgfnodeconnline{Node8}{Node9}
\pgfnodeconnline{Node8}{Node10}
\pgfnodeconnline{Node8}{Node11}
\pgfnodeconnline{Node6}{Node11}
\pgfnodeconnline{Node7}{Node11}
\pgfputat{\pgfxy(6.9,-.4)}{\pgfbox[left,center]{$y_1$}}
\pgfputat{\pgfxy(6, 1.3)}{\pgfbox[left,center]{$y_2$}}
\pgfputat{\pgfxy(5.8, 0.2)}{\pgfbox[left,center]{$y_3$}}
\pgfputat{\pgfxy(5.5, -.6)}{\pgfbox[left,center]{$y_4$}}
\pgfputat{\pgfxy(5.5, -1.2)}{\pgfbox[left,center]{$y_5$}}
\pgfputat{\pgfxy(9.6, 1.3)}{\pgfbox[left,center]{$y_6$}}
\pgfputat{\pgfxy(9.8, -1.3)}{\pgfbox[left,center]{$y_7$}}
\pgfputat{\pgfxy(9, .2)}{\pgfbox[left,center]{$y_8$}}
\pgfputat{\pgfxy(10.1, .6)}{\pgfbox[left,center]{$y_9$}}
\pgfputat{\pgfxy(10.1, -.6)}{\pgfbox[left,center]{$y_{10}$}}
\pgfputat{\pgfxy(12.2, .3)}{\pgfbox[left,center]{$y_{11}$}} 
\end{pgfpicture}
\end{align*} \hskip .9cm \hskip .9cm \hskip .9cm \hskip .9cm \hskip .9cm \hskip .9cm \hskip .3cm Figure 4.  $G_2 =\Gamma(\mathbb{Z}_2 \times \mathbb{Z}_8)$.\\

In Example 2.1, we see that the graph $G_1$ realized by a ring $\mathbb{Z}_{12}$ admits a total perfect code, whereas in Example 2.2, the graph $G_2$ realized by a ring $\mathbb{Z}_2 \times \mathbb{Z}_{8}$ as its zero-divisor graph does not admit a total perfect code. So, it is interesting to characterize rings of which the realized zero-divisor graphs admit total perfect codes.\\

A commutative ring $R$ is called Noetherian if for every ascending chain of ideals $I_1 \subseteq I_2 \subseteq I_3 \subseteq \cdots$, there is a positive integer $m$ such that $I_m = I_{m+k}$ for each positive integer $k$. An Artinian commutative ring is a ring where for every descending chain of ideals $J_1 \supseteq J_2 \supseteq J_3 \supseteq \cdots$, there is an integer $s$ such that $J_s = J_{s+k}$ for each positive integer $k$. Trivially, any finite ring is both Artinian and Noetherian.\\

Let $R$ be a commutative Artinian ring with unity. Then either $R$ is local or $R = R_1  \times \cdots \times R_r \times \mathbb{F}_1 \times \cdots \times \mathbb{F}_s$ , where each $R_i$ is a commutative Artinian local ring with unity that is not a field, each $\mathbb{F}_i$ is a field, and $r$ and $s$ are non-negative integers such that $r + s \geq 2$. The proof of this fact is well known and is a corollary to the Chinese Remainder Theorem. We call this isomorphic image of $R$ the Artinian decomposition of $R$ (we allow $r = 0$ if $s \geq 2$, or $r = 0$ if $s \geq 1)$. Since we need not to consider the case where R is a field, we have three cases to consider: the reduced case (where $r = 0$ and $s \geq 2$ in the Artinian decomposition), the local case (where $n = 1$ and $m = 0$ in the Artinian decomposition), and the mixed cases (not local and not reduced; that is, either $r \geq 1$ and $s \geq 1$, or $r\geq 2$ and $s = 0$ in the Artinian decomposition).\\

Throughout, $(R, m)$ is a finite commutative local ring with unity and with maximal ideal $m \neq {0}$ (that is, $R$ is not a field). Since every element of $R$ is either a unit or a zero-divisor, $Z(R) = m$ (so, in particular, $Z(R)$ is an ideal). Note also that in this case, $|R| = p^n$ for some prime $p$ and an integer $n\geq 2$.\\

A vertex $x$ in a graph $G$ is said to be a \textit{cut-vertex} if the graph resulting by removing the vertex $x$ from $G$ is no longer connected.\\

The following results are characterizations for determining total perfect codes, degree one vertices and cut-vertices in zero-divisor graphs realized by local commutative rings.\\

\begin{proposition} Let $R$ be a finite commutative local ring with unity, then $\Gamma(R)$ admits a total perfect code if and only if $\Gamma(R)$ has degree one vertices.\end{proposition}
{\bf Proof.} Suppose $\Gamma(R)$ admits a total perfect code $C(R)$. Since $R$ is a local ring, therefore by [Corollary 2.7, \cite{AL}] there is some vertex $x\in Z^*(R)$ which is adjacent to all other vertices in $\Gamma(R)$. This implies that $x\in C(R)$. If $\Gamma(R)$ has no degree one vertices, then every vertex of $\Gamma(R)$ has degree more than $1$. In particular, the vertex $y\in C(R)$ distinct from $x$ is adjacent to some other vertex $z\in Z^*(R)$. That is, $yz = 0$ for some $z \neq x$. On the other hand, we have $zx = 0$, since $R$ is local. Thus we have $|Ann(x) \cap C(R)| > 1$, a contradiction. Therefore the graph $\Gamma(R)$ contains atleast one degree one vertex.

Conversely, suppose $\Gamma(R)$ has at least one degree one vertex say $a$. Clearly $a$ is adjacent to $x$, since $x$ is adjacent to all other vertices of $\Gamma(R)$, but $a$ is not adjacent to any other vertex of $\Gamma(R)$. Thus vertices $a$ and $x$ cover all vertices of $\Gamma(R)$. Therefore there is a subset $C(R)$ of $Z^*(R)$  such that $|Ann(x) \cap C(R)| = 1$ for all $x\in Z^*(R)$.\qed

\begin{corollary} Let $R$ be a finite commutative local ring with unity. If $z\in C(R)$ such that $|ann(z)| > 2$ ($z$ is not a degree one vertex in $\Gamma(R)$), then $z$ is a cut-vertex of $\Gamma(R)$.\end{corollary} {\bf Proof.} Proof follows by [Theorem 3, \cite{Rd1}].\qed

\begin{remark} Proposition 2.3 is not true for all finite commutative rings which can be  easily verified from Example 2.2 above. Moreover, the converse of Corollary 2.4 is also not true, since the zero-divisor graphs realized by local commutative rings: $\mathbb{Z}_4[X, Y]/(X^2, Y^2, XY - 2, 2X, 2Y)$, $\mathbb{Z}_2[X, Y]/(X^2, Y^2)$, $\mathbb{Z}_4[X]/(X^2)$, $\mathbb{Z}_4[X]/(X^2 + 2X)$, $\mathbb{Z}_8[X]/(2X, X^2 + 4)$, $\mathbb{Z}_2[X, Y]/(X^2, Y^2 - XY)$, $\mathbb{Z}_4[X, Y]/(X^2, Y^2 - XY, XY - 2, 2X, 2Y)$ has cut-vertices, but does not admit total perfect codes. In fact these are the only local commutative rings which does not admit total perfect codes.
\end{remark}

\begin{proposition} Let $R$ be a finite commutative local ring. Then $\Gamma(R)$ admits a total perfect code if and only if either (a) there is some $x \in R$ such that $|ann(x)| = 2$, where $|Z(R)|\geq 3$ or (b) $R$ is isomorphic to $\mathbb{Z}_9$ or $\mathbb{Z}_3[X]/(X^2)$.\end{proposition}  
{\bf Proof.} Suppose $\Gamma(R)$ admits a total perfect code. By Proposition 2.3, $\Gamma(R)$ has degree one vertices. Let $x$ be a degree one vertex in $\Gamma(R)$. Then either $\Gamma(R)$ has a cut-vertex or $\Gamma(R)$ has only two vertices, implying $R$ is isomorphic to $\mathbb{Z}_9$ or $\mathbb{Z}_3[X]/(X^2)$. 

Let $y$ be a vertex adjacent to $x$. It is clear that $y$ is a cut vertex of $\Gamma(R)$. Therefore $ann(x) = \{0, y\}$ or $ann(x) = \{0, y, x\}$. If $ann(x) = \{0, y, x\}$, then $x^2 = 0$. Note
that $x(x + y) = 0$ and so $x + y \in ann(x)$. The only possibility is that $x + y= 0$. That
is, $x = - y$. However, since $x$ has degree one, the entire graph must consist of only
the vertices $x$ and $y$. This is a contradiction, since a graph must have at least three
vertices to have a cut-vertex. 

The other implication is obvious.\qed

The following result is the immediate consequence of Corollary 2.4, Remark 2.5 and Proposition 2.6.

\begin{corollary} Let $R$ be finite commutative local ring. Then $\Gamma(R)$ has cut-vertices if and only if either (a) there is some $x \in R$ such that $|ann(x)| = 2$ or (b) $R$ is isomorphic to one of the seven rings:  $\mathbb{Z}_4[X, Y]/(X^2, Y^2, XY - 2, 2X, 2Y)$, $\mathbb{Z}_2[X, Y]/(X^2, Y^2)$, $\mathbb{Z}_4[X]/(X^2)$, $\mathbb{Z}_4[X]/(X^2 + 2X)$, $\mathbb{Z}_8[X]/(2X, X^2 + 4)$, $\mathbb{Z}_2[X, Y]/(X^2, Y^2 - XY)$, $\mathbb{Z}_4[X, Y]/(X^2, Y^2 - XY, XY -2, 2X, 2Y)$.\end{corollary}

A ring $R$ is said to be \textit{reduced} if $R$ is free from nilpotent elements. We below investigate total perfect codes in the zero-divisor graph realized by a ring of the form $R = \mathbb{F}_1 \times \cdots \times \mathbb{F}_k$.

\begin{lemma} Let $R$ be a commutative ring with unity. For $k\geq 1$, 

(i) if $R$ realizes a $k$-regular  graph $\Gamma(R)$ on $|Z^*(R)|$ vertices as its zero-divisor graph, then $R$ is a reduced ring,

(ii) if $\Gamma(R)$ admits a total perfect code, then $|Z^*(R)| \equiv 0(mod ~2)$.\end{lemma}

{\bf Proof.} (i) Assume $\Gamma(R)$ is $k$-regular but not a complete graph and assume there is some $x \in Z^*(R)$ such that $x^2 = 0$. Since $\Gamma(R)$ is regular but not complete, there is some $y \in Z^*(R)$ such that $y$ is not adjacent to $x$. Thus, $xy \neq 0$. Also, $xy \neq y$, since $xy \in ann(x)$ but $y\not\in ann(x)$. Clearly, $ann(y) \subseteq ann(xy)$ and $x \in ann(xy)$ but $x \notin ann(y)$. Since $\Gamma(R)$ is $k$-regular, $deg(x) = k$ for all $x \in Z^*(R)$ and therefore $|ann(x)| = k + 2$ if $x^2 = 0$ and $|ann(x)| = k + 1$ otherwise. Thus, $k + 1 \leq |ann(y)| < |ann(xy)| = k + 2$. Therefore $ann(y) = ann(xy) \setminus \{x\}$. It follows that $xy \in ann(xy)$ implies $xy \in ann(y)$. Consequently, this implies $xy \neq x$. So, if $Ann (y) = \{xy, y_1 , y_2, \cdots,  y_{k-1} \}$, then $\{x, y, y_1, y_2, \cdots, y_{k-1} \} = Ann(xy)$, implying $deg(xy) > deg(y)$, a contradiction. Hence, if $\Gamma(R)$ is regular but not complete, then $R$ is reduced.

(ii) Follows because of the facts that $C(R)$ consists of even number of vertices of $\Gamma(R)$ and $k|C(R)| = |Z^*(R)|$.\qed

\begin{lemma} Let $R$ be a commutative ring with unity and let $R$ realizes a complete graph $\Gamma(R)$ on 
 $|Z^*(R)| \geq 2$  vertices as its zero-divisor graph. Then $\Gamma(R)$ admits a total perfect code if and only if $|Z^*(R)| = 2$.
\end{lemma}
{\bf Proof.} For $|Z^*(R)| = 2$, $\Gamma(R)$ is a simple path on two vertices which trivially admit a total perfect. If $|Z^*(R)| > 2$, then we have $|Ann(x) \cap C(R)| > 1$, for every subset $C(R)$ of $Z^*(R)$ and for all $x\in Z^*(R)$. Thus for $|Z^*(R)| \geq 3$, $\Gamma(R)$  does not admit a total perfect code.\qed

\begin{lemma} Let $R$ be a commutative ring with unity and let $R$ realizes a bipartite graph $\Gamma(R)$ as its zero-divisor graph. Then $\Gamma(R)$ admits a total perfect code if and only if it is complete.
\end{lemma}
{\bf Proof.} Let $\Gamma(R)$ be a complete bipartite graph realized by $R$ of order $n + m$, where $n, m \in \mathbb{N}$ with bipartition of $Z^*(R)$ as  $Z^*_{1} (R) = \{x_1 , x_2 , \cdots, x_n \}$ and $Z^*_{2}(R) = \{y_1 , y_2 , \cdots , y_m \}$. For any two vertices  $x_i \in Z^*_{1} (R)$  and $y_j\in Z^*_{2}(R)$ with $1 \leq i \leq n$ and $1 \leq j \leq m$, set $C(R) = \{x_i , y_j \}$. It is easy to check that for
all $a \in Z^*_{1} (R)$ and $b\in Z^*_{2}(R)$, we have $|Ann(a) \cap C(R)| = 1$ and $|Ann(b) \cap C(R)| = 1$. Therefore $C(R)$ is a total perfect code for $\Gamma(R)$.

Suppose that a bipartite graph $\Gamma(R)$ realized by $R$ with partite sets of $Z^*(R)$ as $Z^*_{1}(R)$ and $Z^*_{2}(R)$ admits a total perfect code $C(R)$. It is clear that $C(R)$ contains one vertex from $Z^*_{1}(R)$ and another from $Z^*_{2}(R)$, otherwise $x_px_q = 0$ for some $x_p, x_q\in Z^*_{1}(R)$ or $y_ry_s = 0$ for some $y_r, y_s\in Z^*_{2}(R)$, that is, we would have adjacencies of vertices in a partite set $Z^*_{1}(R)$ or $Z^*_{2}(R)$, which is not possible. If $\Gamma(R)$ is not complete, then there is some vertex $x \in Z^*(R)$ in some partite set which is not adjacent to all vertices of other partite set. If $x \in Z^*_1(R)$, then $x = x_i$ for some $x_i \in Z^*_{1} (R)$. By our assumption $x_i$ is not adjacent to all vertices $y_j \in Z^*_1(R)$, which implies $|Ann(x_i ) \cap C(R)| = \emptyset$, a contradiction. Similarly if $x = y_j$ for some $y_j \in Z^*_2(R)$, then $|Ann(y_j) \cap C(R)| = \emptyset$, again a contradiction to our supposition that $C(R)$ is a total perfect code of $\Gamma(R)$. Thus we conclude that $\Gamma(R)$ is a complete bipartite graph.\qed

Note that Lemma 2.10 is true in general for all bipartite graphs. So, among all bipartite graphs only complete bipartite graphs admit total perfect codes.\\

Let $R$ be a finite commutative reduced ring with unity. If $k$ is the smallest positive integer such that $|R| < 2^k$, then $R$ is a product of $k - 1$ or fewer fields. In fact the smallest ring that is the product of $k$ fields is the finite Boolean ring $\prod\limits_{i=1}^k \mathbb{Z}_2$, which has $2^k$ elements and $2^k - 2$ nonzero zero-divisors. Therefore, it follows that if $R$ has $n$ nonzero zero-divisors and $k$ is the smallest positive integer such that $n < 2^k - 2$, then $R$ is a product of $k - 1$ or fewer fields. Moreover, by [Proposition 2.2, \cite{Rd2}] it is known for a commutative reduced ring $R$ with unity that if $R$ has $k$ maximal ideals, then $R$ is a product of $k$ fields.\\

\begin{proposition} Let $R  = \mathbb{F}_1 \times \cdots \times \mathbb{F}_k$ be a finite commutative reduced ring. Then $\Gamma(R)$ admits a total perfect code if and only if $k = 2$.\end{proposition}

{\bf Proof.} Suppose $k = 2$. If $\mathbb{F}_1 = \mathbb{F}_2 = \mathbb{Z}_2$, then $\Gamma(\mathbb{Z}_2 \times \mathbb{Z}_2)$ is a path on two vertices, which trivially admit a total perfect code. In all other cases, $\Gamma(\mathbb{F}_1 \times \mathbb{F}_2 )$ is either a star graph or a complete bipartite graph with $|\mathbb{F}_1 | + |\mathbb{F}_2 | - 2$ vertices. Therefore by Lemma 2.10, $\Gamma(\mathbb{F}_1 \times \mathbb{F}_2)$ admits a total perfect code.

Conversely, suppose $\Gamma(R)$ admits $C(R)$ as a total perfect code and let $W = \{(u_1, 0, \cdots, 0), (0,$

\noindent $u_2, ,\cdots, 0) ,\cdots,(0, 0, \cdots, u_k)\}$ be a maximal clique in $\Gamma(R)$, where $u_i\in \mathbb{F}_i$, $1\leq i \leq k$. By Lemma 2.9, there is no element of $W$ which is a member of $C(R)$. Moreover, for $k\geq 3$, every vertex $(u_i, 0, \cdots , 0)$ of $W$ has an  adjacent vertex $(0, u_2, \cdots, u_k)$ of degree exactly one (for example, $(1, 0, \cdots ,0)$ would be adjacent to $(0, 1, 1, \cdots, 1)$). Thus there is no subset $C(R)$ of $Z^*(R)$ which covers all vertices of $\Gamma(R)$, that is, there is no subset $C(R)$ of $Z^*(R)$ such that $|Ann(x) \cap Z^*(R)| = 1$ for all $x\in Z^*(R)$, a contradiction to our supposition. Hence, we conclude that for $k\geq 3$, $\Gamma(R)$ does not admit a total perfect code.\qed

\begin{corollary} For $k\geq 3$, the zero-divisor graph $\Gamma(\prod\limits_{i=1}^k \mathbb{Z}_2)$ realized by a finite Boolean ring $\prod\limits_{i=1}^k \mathbb{Z}_2$ does not admit a total perfect code.\end{corollary} 
 
In the remaining of this paper, we consider the rings which are not local and not reduced. To illustrate  this case, given below are some forms of rings, along with the count of non-zero zero-divisors of that structure and an example of the smallest commutative ring with unity having such a structure.\\ 

 $R_1 \times \mathbb{F}$, where $|Z^*(R_1\times \mathbb{F})| = |R_{1}^*| + |\mathbb{F}^*| + |Z^*(R_1)||\mathbb{F}^*|$. The zero-divisor graph realized by a smallest commutative ring of this form is $\Gamma(\mathbb{Z}_4 \times \mathbb{Z}_2)$ with $5$ number of vertices.\

 $R_1 \times R_2$, where $|Z^*(R_1 \times R_2)| = |R_1^{*}| + |R_2^{*}| + |Z^*(R_1)||Z^*(R_2)|$. The zero-divisor graph realized by a smallest commutative ring of this form is $\Gamma(\mathbb{Z}_4 \times \mathbb{Z}_4 \times \mathbb{Z}_4)$ with $11$ number of vertices.\

 $R_1 \times \mathbb{F}_1 \times \mathbb{F}_2$, where $|Z^*(R_1 \times \mathbb{F}_1 \times \mathbb{F}_2)| = |R_{1}^*| + |\mathbb{F}_{1}^*| + |\mathbb{F}_{2}^*| + |R_{1}^*||\mathbb{F}_{1}^*| + |R_{1}^*||\mathbb{F}_{2}^*| + |\mathbb{F}_{1}^*||\mathbb{F}_{2}^*| + |Z^*(R_1)||\mathbb{F}_{1}^*||\mathbb{F}_{2}^*|$. The zero-divisor graph realized by a smallest commutative ring of this form is $\Gamma(\mathbb{Z}_4 \times \mathbb{Z}_2 \times \mathbb{Z}_2)$ with $13$ number of vertices.\

 $R_1 \times {R}_2 \times \mathbb{F}$, where $|Z^*(R_1 \times {R}_2 \times \mathbb{F})| = |R_{1}^*| + |R_{2}^*| + |\mathbb{F}| + |R_{1}^*||{R}_{2}^*| + |R_{1}^*||\mathbb{F}^*| + |R_{2}^*||\mathbb{F}^*| + |Z^*(R_1)||\mathbb{F}^*||{R}_{2}^*| + |Z^*(R_2)||\mathbb{F}^*||{R}_{1}^*| -  |Z^*(R_1)||Z^*(R_2)||\mathbb{F}^*|$. The zero-divisor graph realized by a smallest commutative ring of this form is $\Gamma(\mathbb{Z}_4 \times \mathbb{Z}_4 \times \mathbb{Z}_2)$ with $27$ number of vertices.
 
$R_1 \times \mathbb{F}_1 \times \mathbb{F}_2 \times \mathbb{F}_3$, where $|Z^*(R_1 \times \mathbb{F}_1 \times \mathbb{F}_2 \times \mathbb{F}_3)| = |R_{1}^*| + |\mathbb{F}_{1}^*| + |\mathbb{F}_{2}^*| + |\mathbb{F}_{3}^*| + |R_{1}^*||\mathbb{F}_{1}^*| + |R_{1}^*||\mathbb{F}_{2}^*| + |R_{1}^*||\mathbb{F}_{3}^*| + |\mathbb{F}_{1}^*||\mathbb{F}_{2}^*| + |\mathbb{F}_{1}^*||\mathbb{F}_{3}^*| + |\mathbb{F}_{2}^*||\mathbb{F}_{3}^*| + |R_{1}^*||\mathbb{F}_{1}^*||\mathbb{F}_{2}^*| + |R_{1}^*||\mathbb{F}_{1}^*||\mathbb{F}_{3}^*| + |R_{1}^*||\mathbb{F}_{2}^*||\mathbb{F}_{3}^*| + |\mathbb{F}_1^*||\mathbb{F}_{2}^*||\mathbb{F}_{3}^*| + |Z^*(R_1)||\mathbb{F}_{1}^*||\mathbb{F}_{2}^*||\mathbb{F}_{3}^*|$.The zero-divisor graph realized by a smallest commutative ring of this form is $\Gamma(\mathbb{Z}_4 \times \mathbb{Z}_2 \times \mathbb{Z}_2 \times \mathbb{Z}_2)$ with $29$ number of vertices.

$R_1 \times {R}_2 \times R_3$, where $|Z^*(R_1 \times R_2 \times R_3)| = |R_{1}^*| + |R_{2}^*| + |R_3^*| + |R_{1}^*||{R}_{2}^*| + |R_{1}^*||R_3^*| + |R_{2}^*||R_3^*| + |Z^*(R_1)||R_1^*||{R}_{2}^*| + |Z^*(R_2)||R_3^*||{R}_{1}^*| + |Z^*(R_3)||R_3^*||{R}_{2}^*| -  (|Z^*(R_1)||Z^*(R_2)||R_3^*| + |Z^*(R_1)||Z^*(R_3)||R_2^*| + |Z^*(R_2)||Z^*(R_3)||R_1^*| |Z^*(R_1)||Z^*(R_2)||Z^*(R_3)|)$. The zero-divisor graph realized by a smallest commutative ring of this form is $\Gamma(\mathbb{Z}_4 \times \mathbb{Z}_4 \times \mathbb{Z}_4)$ with $59$ number of vertices.\\

These structures represents only some of the rings, similarly we have other forms of rings in this case. Note that $\mathbb{F}_i$ denotes a finite field and $R_i$ denotes a commutative local ring with unity that is not a field. Moreover, these formulas and any similar formulas for larger rings, are somewhat inductive. That is, one needs to know the structure of all smaller local rings in order to construct the mixed case (not local and not reduced) rings of larger size.\\

In the following result, we discuss total perfect codes in zero-divisor graphs realized by any commutative ring $R$ with unity $1$.\\
 
\begin{theorem} Let $R \cong R_1 \times R_2 \times \cdots \times R_m \times \mathbb{F}_1 \times \mathbb{F}_2 \times \cdots \times \mathbb{F}_n$ be any commutative ring with unity, where each $R_i$ is a finite commutative local ring and each $\mathbb{F}_i$ is a finite field. If $m + n \leq 2$, for $m, n\in \mathbb{Z}_{\geq 0}$ with $m \neq 2$, then $\Gamma(R)$ admits a total perfect code.\end{theorem}

{\bf Proof.} If $m = n = 0$, then there is nothing to prove. Let $n = 0$ and $m = 1$. Then $R$ is a local ring and the result follows from Proposition 2.3. If $m = 0$ and $n \geq 1$, then $R$ is a reduced ring and the result follows from Proposition 2.11.

For $m \geq 1$ and $n \geq 1$, we consider the following cases;

Case 1. $m = 1, n \geq 1$. Then $R \cong R_1 \times \mathbb{F}_1 \times \mathbb{F}_2 \times \cdots \times \mathbb{F}_n$. It is clear that for any unit $u_1 \in R_1$, $W = \{(u_1, 0, \cdots, 0), (0, u_2,\cdots, 0), \cdots, (0, 0, \cdots, u_n)\}$ is a maximal clique in $\Gamma(R)$, where $u_i\in \mathbb{F}_i$, $2\leq i \leq n$, and also every vertex of $W$ has an adjacent vertex $(0, u_2, \cdots, u_k)$ of degree exactly one. Thus for all $x\in Z^*(R)$ there is no subset $C(R)$ of $Z^*(R)$ such that $|C(R) \cap Z^*(R)| = 1$, that is, there is no subset $C(R)$ of $Z^*(R)$ which cover all vertices of $\Gamma(R)$.

Case 2. $m = 1, n = 1$. Then $R \cong R_1 \times \mathbb{F}$. Suppose $|R_1| = k$ and  $|\mathbb{F}| = q$. Let $0, 1, x_{1}, x_{2}, \cdots, x_{k-2}$ be $k$ elements of $R_1$ and $0, 1, z_{1}, z_{2}, \cdots, z_{q - 2}$ be $q$ elements $\mathbb{F}$. Since $R_1$ is a local ring, therefore by [Corollary 2.7, \cite{AL}],  there is some vertex $x_r\in Z^*(R_1)$ which is adjacent to all vertices in $\Gamma(R_1)$. Therefore for all $x_{s_i} \in Z^*(R_1)$ with $1 \leq i \leq |Z^*(R_1)| - 1$, the vertex $(x_r, 0)$ in $\Gamma(R)$ is adjacent to all vertices of the form $(x_{s_i}, 0)$,$(x_{s_i}, 1)$, $(x_{s_i}, z_1)$, $\cdots$, $(x_{s_i}, z_{q-2})$, where $x_r \neq x_{s_i}$ and to all vertices of the form $(0, 1), (0, z_1), \cdots (0, z_{q-2})$. Furthermore, vertices $(0, 1), (0, z_1), \cdots (0, z_{q-2})$ are adjacent to the vertices of the form $(x_{s_i}, 0)$. For $|Z^*(R_1)| \geq 3$, if $\Gamma(R)$ admits $C(R)$ as a total perfect code, then clearly vertices $(x_r, 0)$ and $(0, z)$ are the members of $C(R)$, which implies that $|Ann(x_{s_i}, 0) \cap C(R)| > 1$, a contradiction. Thus $\Gamma(R)$ does not admit a total perfect code. If $|Z^*(R_1)| \leq 2$, then for $z\in \mathbb{F}\setminus \{0\}$ and $x^*_{s_i}\in Z^*(R_1)$, it can be easily verified that there is at least one vertex of the form $(x^*_{s_i}, z)$, which is adjacent to $(x_r, 0)$ but not to any other vertex of $\Gamma(R)$. We set $C(R) = \{(x_r, 0), (x^*_{s_i}, z)\}$ and therefore for all $x\in Z^*(R)$, we have $|Ann(x) \cap C(R)| = 1$. Hence, we conclude that $\Gamma(R)$ admits a total perfect code if and only if $|Z^*(R_1)|\leq 2$.

Case 3. $m = 2, n = 0$. Then $R \cong R_1 \times R_2$. Suppose $|R_1| = |R_2| = k$ and let $0, 1, x_{1}, x_{2}, \cdots, x_{k-2}$ be $k$ elements of $R_1$ and let $0, 1, y_{1}, y_{2}, \cdots, y_{k-2}$ be $k$ elements of $R_2$. We list units in $R_1$ as $1, x_{1}, x_{2}, \cdots, x_{l}$ and non-zero zero-divisors as $x_{l + 1}, x_{l + 2}, \cdots, x_{k-2}$. Also, let $1, y_{1}, y_{2}, \cdots, y_{l}$ be units and $y_{l + 1}, y_{l + 2}, \cdots, y_{k-2}$ be non-zero zero-divisors of $R_2$. Since $R_1$ is local, so there is an element $x_t \in Z^*(R_1)$ with $l+1\leq t\leq k-2$, such that $x_tx_h = 0$ for all $x_h\in Z^*(R_1)$ with $l+1 \leq h \leq k - 3$. Therefore for all $x_h\in Z^*(R_1)$, the vertex $(x_t, 0)$ is adjacent to all vertices of the form $(x_h, 0)$, $(x_h, 1)$, $(x_h, x_1)$, $\cdots$, $(x_h, x_{k-2})$ in $\Gamma(R)$. Similarly for all $y_h\in Z^*(R_2)$ with $l+1 \leq h \leq k - 3$, there is an element $y_t\in Z^*(R_2)$ such that the vertex $(0, y_t)$ is adjacent to all vertices of the form $(0, y_{h})$, $(1, y_h)$, $(y_1, y_h)$ $\cdots$, $(y_{k-2}, y_h)$ in $\Gamma(R)$. Furthermore, for all $y\in R_2 \setminus\{0\}$, every vertex of the form $(0, y)$ is adjacent to the vertex $(x_t, 0)$ and for all $x\in R_1 \setminus\{0\}$, every vertex of the form $(x, 0)$ is adjacent to $(0, y_t)$. If $\Gamma(R)$ admits $C(R)$ as a total perfect code, then $(x_t, 0)$ and $(0, y_t)$ are the members of $C(R)$, since $(x_t, 0)$ and $(0, y_t)$ cover all vertices of $\Gamma(R)$, but for all $y_h\in Z^*(R_2)$, we have $|Ann(0, y_h) \cap C(R)| > 1$, a contradiction. Thus $\Gamma(R)$ does not admit a total perfect code.

Case 4. $m = 2, n\geq 1$. Then $R \cong R_1 \times R_2 \times \mathbb{F}_1 \times \mathbb{F}_2 \times \cdots \times \mathbb{F}_n$. By cases 1 and 3, there is no subset $C(R)$ of $Z^*(R)$ such that for all $x\in R$, $|Ann(x) \cap C(R)| = 1$. 

Case 5. $m \geq 3, n = 0$. Then $R \cong R_1 \times R_2 \times \cdots \times R_m$. Again by case 3, we see that there is no subset $C(R)$ of $Z^*(R)$ such that for all $x\in R$, $|Ann(x) \cap C(R)| = 1$. Therefore $\Gamma(R)$ does not admit a total perfect code.\qed

We present below some examples of graphs realized as zero-divisor graphs which does not admit total perfect codes. These examples illustrate cases $2$, $3$ and $4$ of the preceding theorem.\\

\begin{example} Consider a ring $\mathbb{Z}_2 \times \mathbb{Z}_2[X, Y]/(X, Y)^2$, where $\mathbb{Z}_2$ is a reduced ring and $\mathbb{Z}_2[X, Y]/(X, Y)^2$ is a local ring. It is clear by the structure of a ring that $|Z^*(\mathbb{Z}_2[X, Y] /(X, Y)^2)| > 2$. The zero-divisor graph realized by a ring $\mathbb{Z}_2[X, Y]/(X, Y)^2$ is shown in figure 3. It can be easily verified from the graph that $\Gamma(\mathbb{Z}_2 \times \mathbb{Z}_2[X, Y]/(X, Y)^2)$ does not admit a total perfect code.\end{example}
\begin{align*}
\begin{pgfpicture}{9.5cm}{-1cm}{2cm}{.5cm}
\pgfnodecircle{Node1}[fill]{\pgfxy(7,0)}{0.1cm}
\pgfnodecircle{Node2}[fill]{\pgfxy(6,1)}{0.1cm}
\pgfnodecircle{Node3}[fill]{\pgfxy(6, .5)}{0.1cm}
\pgfnodecircle{Node4}[fill]{\pgfxy(6, -.5)}{0.1cm}
\pgfnodecircle{Node5}[fill]{\pgfxy(6, -1)}{0.1cm}
\pgfnodecircle{Node6}[fill]{\pgfxy(11, 1)}{0.1cm}
\pgfnodecircle{Node7}[fill]{\pgfxy(11,-1)}{0.1cm}
\pgfnodecircle{Node8}[fill]{\pgfxy(10, 0)}{0.1cm}
\pgfnodecircle{Node9}[fill]{\pgfxy(13, 0)}{0.1cm}
\pgfnodecircle{Node10}[fill]{\pgfxy(13,1)}{0.1cm}
\pgfnodecircle{Node11}[fill]{\pgfxy(13,-1)}{0.1cm}
\pgfnodeconnline{Node1}{Node2}
\pgfnodeconnline{Node1}{Node3}
\pgfnodeconnline{Node1}{Node4}
\pgfnodeconnline{Node1}{Node5}
\pgfnodeconnline{Node1}{Node6}
\pgfnodeconnline{Node1}{Node7}
\pgfnodeconnline{Node1}{Node8}
\pgfnodeconnline{Node6}{Node8}
\pgfnodeconnline{Node6}{Node9}
\pgfnodeconnline{Node6}{Node7}
\pgfnodeconnline{Node8}{Node7}
\pgfnodeconnline{Node8}{Node9}
\pgfnodeconnline{Node8}{Node10}
\pgfnodeconnline{Node8}{Node11}
\pgfnodeconnline{Node6}{Node10}
\pgfnodeconnline{Node6}{Node11}
\pgfnodeconnline{Node7}{Node10}
\pgfnodeconnline{Node7}{Node11}
\pgfnodeconnline{Node7}{Node9}
\pgfputat{\pgfxy(6.9,-.4)}{\pgfbox[left,center]{}}
\pgfputat{\pgfxy(6, 1.1)}{\pgfbox[left,center]{}}
\pgfputat{\pgfxy(6, 2.1)}{\pgfbox[left,center]{}}
\pgfputat{\pgfxy(3.9, -.6)}{\pgfbox[left,center]{}}
\pgfputat{\pgfxy(4.5, -2)}{\pgfbox[left,center]{}}
\pgfputat{\pgfxy(9.6,2.8)}{\pgfbox[left,center]{}}
\pgfputat{\pgfxy(9.8, -2.5)}{\pgfbox[left,center]{}}
\pgfputat{\pgfxy(9.5, .2)}{\pgfbox[left,center]{}}
\pgfputat{\pgfxy(11.8, .9)}{\pgfbox[left,center]{}}
\pgfputat{\pgfxy(11.8, -.3)}{\pgfbox[left,center]{}}
\pgfputat{\pgfxy(14.1, .2)}{\pgfbox[left,center]{}} 
\end{pgfpicture}
\end{align*} \hskip .9cm \hskip .9cm \hskip .9cm \hskip .9cm \hskip .9cm \hskip .9cm $Figure ~3$. $\Gamma(\mathbb{Z}_2 \times \mathbb{Z}_2[X, Y]/(X, Y)^2)$.\

\begin{example} Consider the product $\mathbb{Z}_4 \times \mathbb{Z}_4$ of two commutative local rings. The zero-divisor graph realized by a ring $\mathbb{Z}_4 \times \mathbb{Z}_4$ is shown in figure $4$ below. It is clear from the graph that $\Gamma(\mathbb{Z}_4 \times \mathbb{Z}_4)$ deos not admit a total perfect code.\end{example}
\begin{align*}
\begin{pgfpicture}{9cm}{-2cm}{4cm}{.5cm}
\pgfnodecircle{Node1}[fill]{\pgfxy(6,0)}{0.1cm}
\pgfnodecircle{Node2}[fill]{\pgfxy(6,-1.5)}{0.1cm}
\pgfnodecircle{Node3}[fill]{\pgfxy(7,0)}{0.1cm}
\pgfnodecircle{Node4}[fill]{\pgfxy(7,-1.5)}{0.1cm}
\pgfnodecircle{Node5}[fill]{\pgfxy(8,0)}{0.1cm}
\pgfnodecircle{Node6}[fill]{\pgfxy(8,-1.5)}{0.1cm}
\pgfnodecircle{Node7}[fill]{\pgfxy(9,-.8)}{0.1cm}
\pgfnodecircle{Node8}[fill]{\pgfxy(9, -2.1)}{0.1cm}
\pgfnodecircle{Node9}[fill]{\pgfxy(11, -1.5)}{0.1cm}
\pgfnodecircle{Node10}[fill]{\pgfxy(11,0)}{0.1cm}
\pgfnodecircle{Node11}[fill]{\pgfxy(9,.7)}{0.1cm}
\pgfnodeconnline{Node1}{Node2}
\pgfnodeconnline{Node2}{Node3}
\pgfnodeconnline{Node3}{Node4}
\pgfnodeconnline{Node4}{Node5}
\pgfnodeconnline{Node5}{Node6}
\pgfnodeconnline{Node1}{Node6}
\pgfnodeconnline{Node3}{Node6}
\pgfnodeconnline{Node2}{Node3}
\pgfnodeconnline{Node4}{Node1}
\pgfnodeconnline{Node6}{Node7}
\pgfnodeconnline{Node8}{Node6}
\pgfnodeconnline{Node6}{Node9}
\pgfnodeconnline{Node5}{Node11}
\pgfnodeconnline{Node5}{Node10}
\pgfnodeconnline{Node5}{Node7}
\pgfnodeconnline{Node5}{Node2}
\pgfputat{\pgfxy(6.9,-.4)}{\pgfbox[left,center]{}}
\pgfputat{\pgfxy(6, 1.1)}{\pgfbox[left,center]{}}
\pgfputat{\pgfxy(6, 2.1)}{\pgfbox[left,center]{}}
\pgfputat{\pgfxy(3.9, -.6)}{\pgfbox[left,center]{}}
\pgfputat{\pgfxy(4.5, -2)}{\pgfbox[left,center]{}}
\pgfputat{\pgfxy(9.6,2.8)}{\pgfbox[left,center]{}}
\pgfputat{\pgfxy(9.8, -2.5)}{\pgfbox[left,center]{}}
\pgfputat{\pgfxy(9.5, .2)}{\pgfbox[left,center]{}}
\pgfputat{\pgfxy(11.8, .9)}{\pgfbox[left,center]{}}
\pgfputat{\pgfxy(11.8, -.3)}{\pgfbox[left,center]{}}
\pgfputat{\pgfxy(14.1, .2)}{\pgfbox[left,center]{}} 
\end{pgfpicture}
\end{align*}

\hskip .9cm \hskip .9cm \hskip .9cm \hskip .9cm \hskip .9cm \hskip .9cm \hskip .3cm $Figure ~4$. $\Gamma(\mathbb{Z}_4 \times \mathbb{Z}_4)$.\

\begin{example} Here we consider the case of two local rings and a field. Let $R = \mathbb{Z}_2 \times \mathbb{Z}_2 \times \mathbb{Z}_2[X]/(X)^2$. It can be easily verified from the zero-divisor graph realized by $R$ (see figure 5) that $\Gamma(R)$ deos not admit a total perfect code.\end{example}

\begin{align*}
\begin{pgfpicture}{8.5cm}{-0.1cm}{3cm}{3cm}
\pgfnodecircle{Node1}[fill]{\pgfxy(7,0)}{0.1cm}
\pgfnodecircle{Node2}[fill]{\pgfxy(5.5,1)}{0.1cm}
\pgfnodecircle{Node3}[fill]{\pgfxy(5.5,-1)}{0.1cm}
\pgfnodecircle{Node4}[fill]{\pgfxy(7,-1.5)}{0.1cm}
\pgfnodecircle{Node5}[fill]{\pgfxy(9, -1.5)}{0.1cm}
\pgfnodecircle{Node6}[fill]{\pgfxy(9, 0)}{0.1cm}
\pgfnodecircle{Node7}[fill]{\pgfxy(7,1.5)}{0.1cm}
\pgfnodecircle{Node8}[fill]{\pgfxy(6, 3)}{0.1cm}
\pgfnodecircle{Node9}[fill]{\pgfxy(8,3.5)}{0.1cm}
\pgfnodecircle{Node10}[fill]{\pgfxy(10,2)}{0.1cm}
\pgfnodecircle{Node11}[fill]{\pgfxy(10,1)}{0.1cm}
\pgfnodecircle{Node12}[fill]{\pgfxy(10,-1.5)}{0.1cm}
\pgfnodeconnline{Node1}{Node2}
\pgfnodeconnline{Node1}{Node3}
\pgfnodeconnline{Node1}{Node4}
\pgfnodeconnline{Node1}{Node4}
\pgfnodeconnline{Node1}{Node7}
\pgfnodeconnline{Node1}{Node7}
\pgfnodeconnline{Node1}{Node6}
\pgfnodeconnline{Node1}{Node5}
\pgfnodeconnline{Node11}{Node9}
\pgfnodeconnline{Node1}{Node10}
\pgfnodeconnline{Node4}{Node5}
\pgfnodeconnline{Node4}{Node6}
\pgfnodeconnline{Node5}{Node6}
\pgfnodeconnline{Node8}{Node7}
\pgfnodeconnline{Node7}{Node9}
\pgfnodeconnline{Node1}{Node9}
\pgfnodeconnline{Node7}{Node6}
\pgfnodeconnline{Node7}{Node10}
\pgfnodeconnline{Node6}{Node11}
\pgfnodeconnline{Node11}{Node10}
\pgfnodeconnline{Node6}{Node12}
\pgfnodeconnline{Node11}{Node9}
\pgfnodeconnline{Node7}{Node6}
\pgfputat{\pgfxy(6.9,-.4)}{\pgfbox[left,center]{}}
\pgfputat{\pgfxy(6, 1.1)}{\pgfbox[left,center]{}}
\pgfputat{\pgfxy(6, 2.1)}{\pgfbox[left,center]{}}
\pgfputat{\pgfxy(3.9, -.6)}{\pgfbox[left,center]{}}
\pgfputat{\pgfxy(4.5, -2)}{\pgfbox[left,center]{}}
\pgfputat{\pgfxy(9.6,2.8)}{\pgfbox[left,center]{}}
\pgfputat{\pgfxy(9.8, -2.5)}{\pgfbox[left,center]{}}
\pgfputat{\pgfxy(9.5, .2)}{\pgfbox[left,center]{}}
\pgfputat{\pgfxy(11.8, .9)}{\pgfbox[left,center]{}}
\pgfputat{\pgfxy(11.8, -.3)}{\pgfbox[left,center]{}}
\pgfputat{\pgfxy(14.1, .2)}{\pgfbox[left,center]{}}
\pgfputat{\pgfxy(14.1, -.2)}{\pgfbox[left,center]{}}
\end{pgfpicture}
\end{align*}\\\\\\

\hskip .9cm \hskip .9cm \hskip .9cm \hskip .9cm \hskip .7cm  $Figure ~5$.  $\Gamma(\mathbb{Z}_2 \times \mathbb{Z}_2 \times \mathbb{Z}_2[X]/(X)^2)$.\\

We conclude this paper with the following result in which we determine the cardinality of a total perfect code in $\Gamma(R)$. In fact for every commutative ring $R$ with unity which realizes some $\Gamma(R)$ as its zero-divisor graph, we exhibit that either $C(R) = \phi$ or $C(R)$ is a trivial matching.\\

\begin{theorem} Let $R$ be a commutative ring with unity. If $\Gamma(R)$ admits $C(R)$ as a total perfect code, then $|C(R)| = 2k$, where $k\in \{0, 1\}$.\end{theorem}

{\bf Proof.} By Example 2.2, we see that $\Gamma(R)$ does not always admits a total perfect code. Similarly we have other examples where $C(R) = \phi$. This implies $|C(R)| = 0$. Suppose $\Gamma(R)$ admits a total perfect code. Assume that $|C(R)| > 2$. By definition, $C(R)$ is matching which implies that $C(R)$ consists of even number of vertices. If possible, suppose $|C(R)| = 4$ and let $X = (a_1,a_2)$, $Y = (c_1,c_2)$ be two elements of $C(R)$. It is clear that $X \cap Y = \phi$ and there is an edge $Z = (b_1, b_2)$, where $b_1\in Ann(a_2)$ and $b_2\in Ann(c_1)$ with $b_1 \neq b_2$. Therefore, $ d(a_1, c_1) = d(a_1, b_1) + d(b_1, c_1) = 2 + 2 = 4$,
which is a contradiction, since diameter of $\Gamma(R)$ is atmost 3. Hence, we conclude that if $\Gamma(R)$ admits a total perfect code, then $C(R)$ is a trivial matching and therefore $|C(R)| = 2k$, where $k\in \{0, 1\}$.\qed
\begin{remark} Apart from determining a commutative ring for the given graph, the study of total perfect codes in graphs realized as zero-divisor graphs leads us to one more significance. Using Theorem 2.13, we can construct a family of diameter $3$ graphs realized as zero-divisor graphs which does not admit total perfect codes. $\Gamma(\mathbb{Z}_2 \times \mathbb{Z}_2 \times \mathbb{Z}_2[X]/(X)^2)$, $\Gamma(\mathbb{Z}_4 \times \mathbb{Z}_4)$, $\Gamma(\mathbb{Z}_2 \times \mathbb{Z}_2[X, Y]/(X, Y)^2)$ and $\Gamma(\mathbb{Z}_2 \times \mathbb{Z}_8)$ are some members of this family. J. Kratcohvil in (1986) \cite{K}, proved the non existence of non-trivial perfect codes over complete bipartite graphs. Thus, we have another family of graphs which does not admit total perfect codes. In fact the whole family of graphs realized as zero-divisor graphs does not admit non-trivial total perfect codes.\vskip .7cm
\end{remark}

{\bf Acknowledgements}. I would like to thank the Department of Atomic Energy, Government of India for providing me the financial support under Grant No. HRI/4042/3784 and Harish-Chandra Research Institute, Allahabad for research facilities. Moreover, I would like to thank the referee for his/her constructive comments which helped improving the quality of the paper.

\end{document}